\documentclass[a4paper]{article}

\usepackage{MyLaTeX}
\usepackage[vcentermath,enableskew]{youngtab}
\title{Symmetric Group Character Degrees and Hook Numbers}
\author{David A.~Craven, University of Oxford}
\date{June 2006}

\numberwithin{defn}{section}

\begin{document}
\maketitle

\noindent In this article we prove the following result: that for any two natural numbers $k$ and $\ell$, and for all sufficiently large symmetric groups $S_n$, there are $k$ disjoint sets of $\ell$ irreducible characters of $S_n$, such that each set consists of characters with the same degree, and distinct sets have different degrees. In particular, this resolves a conjecture most recently made by Moret\'o in \cite{moreto}. The methods employed here are based upon the duality between irreducible characters of the symmetric groups and the partitions to which they correspond. Consequently, the paper is combinatorial in nature.

\section{Introduction}

This article will discuss the degrees of irreducible characters of the symmetric group, and will in particular prove the following theorem, utilizing a combinatorial approach.

\begin{thm}\label{mainthm1} Let $k$ and $\ell$ be natural numbers, and let $S_n$ denote the symmetric group on $n$ letters. There exists an integer $N$ such that for all $n\geq N$, there are $k$ disjoint sets of $\ell$ irreducible ordinary characters, each set consisting of characters with the same degree, and distinct sets have different degrees; in other words, as $n$ tends to infinity, the number of disjoint sets of $\ell$ irreducible ordinary characters of $S_n$ all $\ell$ of which have the same degree also tends to infinity.
\end{thm}

Theorem \ref{mainthm1} has the following specialization to the case $k=1$: a conjecture stated, for example, in \cite{moreto}.

\begin{thm}\label{maincor1} Let $m(n)$ denote the size of a largest subset of irreducible characters of $S_n$, each of which has the same degree. Then $m(n)\to \infty$ as $n\to\infty$.
\end{thm}

Brauer's Problem 1 (in \cite{brauer}) asks whether one can determine which $\C$-algebras are isomorphic with group algebras; this is clearly equivalent to asking if, given a multiset $M$ of positive integers, $M$ is the multiset $M_G$ of degrees of irreducible ordinary characters of some finite group $G$. While this question might be too hard in general, a subproblem is to ask whether there is bound on the order of a finite group $G$ in terms of the multiplicities of the irreducible character degrees. Recently, Alexander Moret\'o in \cite{moreto} has shown that the order of any finite group $G$ is bounded by a function of the maximum $m(G)$ of the multiplicities of irreducible character degrees if the corresponding result is true for the symmetric groups only. Firstly, the problem was reduced to the finite simple groups. Discarding the sporadics, this left the groups of Lie type, for which the result is proven by methods of algebraic conjugacy, and the alternating groups, for which the algebraic-conjugacy method obviously fails. Significantly different methods to those employed by Moret\'o are needed to deal with the symmetric groups. Theorem \ref{maincor1} closes this last possible gap, and hence yields the following corollary.

\begin{cor} Let $G$ be a finite group, and let $m(G)$ denote the maximal multiplicity among the degrees of the irreducible ordinary characters. Then $|G|$ is bounded by a function of $m(G)$.
\end{cor}

As was said before, the character theory of the symmetric group can be studied in an entirely combinatorial manner: the Hook Length Formula connects the degrees of characters with products of integers associated with partitions. The rest of this paper will essentially be devoted to a combinatorial proof of Theorem \ref{mainthm1}.

A \emph{partition} of $n$ is a weakly decreasing sequence of positive integers $\lambda=(\lambda_1,\dots,\lambda_r)$ with $\sum \lambda_i=n$. The $\lambda_i$ are called the \emph{parts} of $\lambda$. If some of the $\lambda_i$ are equal, then we will often abbreviate this, so that $(1,1,1)$ becomes $(1^3)$, for example. Recall that irreducible representations of $S_n$ are in one-to-one correspondence with partitions of $n$ (see for instance \cite{jameskerber}); the degree of the representation in correspondence with the partition $\lambda$ is given by the famous formula of Frame, Robinson and Thrall \cite{framerobthrall}
\[ \chi_\lambda(1)=\frac{n!}{\prod_{(i,j)} h(i,j)},\]
where $h(i,j)$ is the \emph{hook number}, which we will define now. Partitions can be represented as tableaux, with the number of boxes in each row equal to the parts of the partition, so that for example $(4,2,2,1)$ is represented as
\[ \yng(4,2,2,1).\]
Then the hook number of a box $a$ of the tableau is simply the sum of the number of boxes below the box $a$, the number of boxes to the right of $a$, and 1 (for the box $a$ itself). Thus the hook numbers of the partition above are
\[ \young(7521,42,31,1),\]
and so the degree of the character corresponding to this partition is
\[ \frac{9!}{1680}=216.\]
Notice that the conjugate partition---the partition reflected in the diagonal running down and to the right---possesses the same hook numbers, and so the corresponding characters have the same degree. (In fact, the `conjugated' character is the tensor product of the original character and the alternating character obtained from the homomorphism $S_n\to \{\pm 1\}$.)

In general, knowing when two multisets of hook numbers have the same product is a difficult problem; however, if two partitions have the same hook number multisets, then they certainly correspond to characters with the same degree. In \cite{hermanchung}, Herman and Chung provide examples of two non-conjugate partitions with equal hook number multisets, namely
\[ A_n=(n+6,n+3,n+3,2)\quad\text{and}\quad A_n'=(n+5,n+5,n+2,1,1),\]
and
\[ B_n=(n+8,n+4,n+3,3,1)\quad\text{and}\quad B_n'=(n+7,n+6,n+2,2,1,1),\]
for all $n\geq 0$. The characters described in Theorem \ref{mainthm1} in fact have the extra condition that the hook numbers of the corresponding partitions are equal as multisets, in addition to their products being equal. It is partitions whose hook numbers are equal as multisets that are the focus here. The statement of the full result that we prove is the following.

\begin{thm}\label{mainthm2} Let $k$ and $\ell$ be natural numbers. Then for all sufficiently large $n$, there are $k$ disjoint sets of $\ell$ partitions of $n$, such that all of the $\ell$ partitions in each set have the same multiset of hook numbers, and distinct sets contain partitions with different hook numbers, and moreover different products of hook numbers.\end{thm}

The proof of this theorem requires a construction of a special type of partition, called the \emph{enveloping partition}. To each partition $\lambda$, we can associate another partition, $E(\lambda)$. The precise construction of $E(\lambda)$ is given in Section \ref{proofofthm}, but the reason these are constructed is that if $\lambda$ and $\mu$ have the same hook numbers, so do $E(\lambda)$ and $E(\mu)$. This construction enables us to build up larger collections of partitions with the same hook numbers from smaller collections, the crux of the proof of Theorem \ref{mainthm2}.

In Theorem \ref{mainthm2}, we said that the distinct sets had different products of hook numbers. To show that they have different multisets of hook numbers is easy, but in order to prove that they have different products, we require the famous Bertrand's Postulate from number theory, namely that if $n$ is any natural number, there exists a prime between $n$ and $2n$.

\section{Definitions and Preliminaries}\label{defns}

For the rest of this article, $\lambda$ and $\mu$ will generally denote partitions. The partition $\lambda$ will normally have $r$ rows and $c$ columns, and we will write $t=r+c$. If $\lambda$ denotes a partition, $\lambda^{(c)}$ will always denote its conjugate. Write $|\lambda|$ for the weight of a partition $\lambda$; i.e., of what number it is a partition. We denote by $H(\lambda)$ the multiset of hook numbers, and by $H_i(\lambda)$ the number of times that $i$ appears in $H(\lambda)$.

There is an equivalence relation $\sim$ on the set of all partitions, where $\lambda\sim\mu$ if and only if $H(\lambda)=H(\mu)$. If $\lambda\sim \mu$, we say that $\lambda$ and $\mu$ are \emph{clustered}, and a collection of partitions $\Lambda$ of $n$ is a \emph{cluster} if it is a subset of an equivalence class of the relation $\sim$. The \emph{size} of a cluster is the cardinality of the cluster; i.e., how many partitions are in the cluster.

We have given two families of clusters of size 4, namely $A_n$, $A_n'$ and their conjugates, and also $B_n$, $B_n'$ and their conjugates. There are infinitely many examples, including
\[ C_n=(n+10,n+4,n+4,4,2)\quad\text{and}\quad C_n'=(n+8,n+8,n+2,2,2,1,1),\]
and two-parameter families such as 
\[ (n+6,(n+4)^r,(n+3)^2,2)\quad\text{and}\quad ((n+5)^2,(n+4)^r,n+2,1^2)\]
for all $n,r\geq 0$. (If $n=r-1$, then the two partitions are actually conjugate, but otherwise they are not.) We see therefore that not all families are parameterized by a single variable; Section \ref{vet} will produce the optimal result in this direction, namely that we have clusters parameterized by all finite strings of non-negative integers. Since the set of all partitions is countable, a countably-parameterized cluster is best-possible.

We shall introduce some terminology to make the discussion of clusters such as $\{A_n,A_n'\}$ given above easier. Consider, for example, $A_4$. This looks like the diagram below.
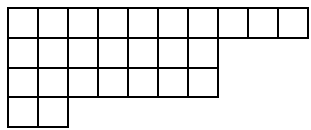
\[\yng(10,7,7,2)\]
Notice that for the partitions in this infinite sequence $(n+6,(n+3)^2,2)$, one section is extended with increasing $n$, and the rest---in this case the partition $(2)$---stays the same. This is entirely typical of the general case that we wish to consider, and so we name some features of these types of sequence.

If $\lambda$ is a partition, and $n$ and $p$ are non-negative integers, denote by $\lambda^{(p,n)}$ the partition obtained from $\lambda$ by incrementing the first $p$ parts by $n$. A cluster $\Lambda$ is called a \emph{periodic} cluster, of \emph{period $p$}, if for all $n\geq 0$, the set
\[ \{\lambda^{(p,n)}:\lambda\in\Lambda\}\]
is also a cluster. For example then, the set $\{A_0,A_0'\}$ is a periodic cluster of period 3. [Note that there exist clusters that are not periodic of any period, but we will not consider such clusters here.]

This definition is, however, not flexible enough for our purposes. For example, suppose that we have a cluster $\Lambda$, consisting of $A_0$, $A_0'$, and their conjugates. The definition of a periodic cluster as it stands would imply that although $\{A_0,A_0'\}$ is periodic, $\Lambda$ is not. We will extend our definition, and say that a cluster $\Lambda$ is periodic if there is a periodic cluster $\Lambda'$ such that for every $\lambda\in \Lambda$, either $\lambda\in \Lambda'$ or $\lambda^{(c)}\in \Lambda'$.

Suppose that $\lambda$ is an element of a period-$p$ cluster. The \emph{rump} of $\lambda$ is the partition consisting of the first $p$ parts of $\lambda$, and the \emph{remainder}, normally denoted $\lambda^{(r)}$, is the partition consisting of all parts of $\lambda$ apart from the first $p$. Finally, the \emph{front section}, normally denoted $\lambda^{(f)}$, is the partition that remains upon deleting all columns of the rump of height $p$ apart from the right-most one. (Then a front section is an arbitrary partition whose smallest part has value 1.) Thus the rump is the piece of the partition that is incremented, the remainder is the piece that does not alter between $\lambda$ and $\lambda^{(p,n)}$, and the front section is the piece of the rump that determines its shape. For example, the rump of $A_4$, pictured above, is $(10,7,7)$, the remainder is $(2)$, and the front section is $(4,1,1)$. Notice that for all $n$, the front sections and remainders of $A_n$ are identical.

We collect all of the information about a partition $\lambda$ in a so-called \emph{$p$-partition datum}. Let $p$ be a natural number less than the number of rows of $\lambda$. The $p$-partition datum consists of three objects: the front section of $\lambda$, the remainder of $\lambda$, and the \emph{extension} $e_\lambda$ of $\lambda$, the quantity $\lambda_p-\lambda_{p+1}$. We write $\lambda=[\lambda^{(f)},\lambda^{(r)},e_\lambda]_p$. In our example of $A_4$,
\[ A_4=[(4,1,1),(2),5]_3.\]
This method of writing partitions makes it much easier to define periodic clusters, which we will have cause to do later.

Now consider another of the period-3 clusters, say $\{B_n,B_n'\}$. Suppose that $\lambda=(n+7,n+6,n+2,2,1^2)$ and $\mu=(n+8,n+4,n+3,3,1)$. Notice firstly that the remainders $\lambda^{(r)}$ and $\mu^{(r)}$ are conjugate. Secondly, notice that if we rotate the partition $\mu^{(f)}$ by $180$ degrees, it fits perfectly against $\lambda^{(f)}$; we will call two front sections that fit together in this way \emph{complementary}. Here is the example of $A_4$.

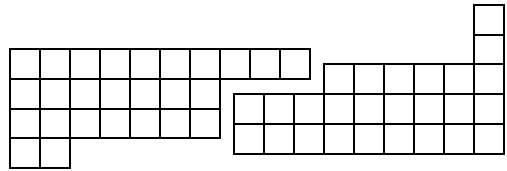
\[ \young(\,\,\,\,\,\,\,\,\,\,,\,\,\,\,\,\,\,,\,\,\,\,\,\,\,,\,\,) \hspace{-0.4in} \raisebox{-0.2in}{\young(::::::::\,,::::::::\,,:::\,\,\,\,\,\,,\,\,\,\,\,\,\,\,\,,\,\,\,\,\,\,\,\,\,)}\]

We have demonstrated two properties of this family of partitions from \cite{hermanchung}. The clusters $\{A_n,A_n'\}$ and $\{C_n,C_n'\}$ also have these properties, as does the two-parameter periodic cluster given above. These two properties are not inherent to a periodic cluster: the periodic clusters
\[ (n+8,(n+5)^2,5,3,2^3)\quad\text{and}\quad((n+7)^2,n+4,4^2,3,1^3),\]
and
\[ (n+8,(n+5)^3,n+3,(n+2)^2,2)\quad\text{and}\quad((n+7)^2,(n+4)^3,n+3,n+1,1^2),\]
are examples of periodic clusters whose remainders are not conjugate, and whose front sections are not complementary, respectively. They do, however, satisfy other conditions weaker than their remainders being conjugate and their front sections being complementary, as we shall see.

Suppose that $\lambda$ and $\mu$ are clustered, of period $p$. Section \ref{splintering}, whose main focus is the Splintering Lemma,  shows that, as long as the two front sections $\lambda^{(f)}$ and $\mu^{(f)}$ are complementary, `most' of the hook numbers in the rumps of $\lambda$ and $\mu$  will automatically be the same. The Extension Lemma is the subject of the succeeding section: this section will give a sufficient condition for a cluster, $\{\lambda,\mu\}$ say, to be periodic. Such a result reduces the task of showing that, for example, $\{C_n,C_n'\}$ is a cluster for all $n$ to showing that $C_0$ and $C_0'$ are clustered, and that they satisfy the hypotheses of the Extension Lemma. The Splintering and Extension Lemmas together constitute the very basic tools used in the proof of Theorem \ref{mainthm2}. A large reduction of the problem, essentially to finding just one cluster of each size for each weight, is the content of the Vertical Expansion Theorem, given in Section \ref{vet}. The proof of Theorem \ref{mainthm2} itself occupies Sections \ref{prepforthm} and \ref{proofofthm}.

If $\lambda=(\lambda_1,\dots,\lambda_r)$ is a partition with $r$ rows and $c$ columns, then the \emph{remnant} $\widetilde{\lambda}$ (the term is specific to this article) is given by
\[ \widetilde{\lambda}=(c-\lambda_r,c-\lambda_{r-1},\dots,c-\lambda_{1}),\]
assuming that the zeros are removed. Pictorially, we can think of the remnant as being the partition that has been removed from the $(r\times c)$-rectangle to create $\lambda$. For example, if $\lambda=(5,5,4,2,1)$, then $\widetilde{\lambda}=(4,3,1)$.
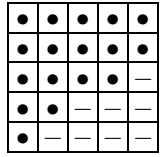
\[ \young(\bullet\bullet\bullet\bullet\bullet,\bullet\bullet\bullet\bullet\bullet,\bullet\bullet\bullet\bullet-,\bullet\bullet---,\bullet----)\]
In the diagram above, the boxes with $\bullet$ in them denote $\lambda$, and those with $-$ in them denote the partition $\widetilde{\lambda}$, rotated by $180$ degrees. Whenever $\lambda$ denotes a partition, $\widetilde{\lambda}$ will always denote its remnant.

For future reference, we give here the basic results about hook numbers that we need. These are well-known, and the reader is invited to give a proof if needed. The last four are very similar in nature, and are included individually here so we can use the specific forms when required.

\begin{lem}\label{oldresultonhooks} Let $\lambda$ be a partition, with $r$ rows and $c$ columns, and write $t=r+c$. Denote by $h(i,j)$ the hook number in the $(i,j)$ box. We have the following:
\begin{enumerate}
\item $h(1,1)=t-1$;
\item $h(i,j)=h(i,1)+h(1,j)-h(1,1)$;
\item we have
\[ \{h(i,1): 1\leq i\leq r\}\cup \{(t-1)-h(1,j): 1\leq j\leq c\}=\{0,1,\dots,t-1\};\]
\item if $A$ denotes the set of hook numbers in the left-hand column, $A'$ denotes the complement of $A$ in $\{0,\dots,t-1\}$, and $B$ denotes the set of hook numbers in the top row, then
\[ B=\{ (t-1)-a': a'\in A'\};\]
\item if $A$ is the set of hook numbers in the left-hand column, $a$ is the left-most hook number in row $i$, and $A'$ is the complement of $A$ in $\{0,\dots,t-1\}$, then the hook numbers in row $i$ are
\[ \{a-a': a'\in A',\; a>a'\};\]
and
\item if $A$ is the set of hook numbers in the left-hand column, and $A'$ denotes its complement in $\{0,\dots,t-1\}$, then
\[ H(\lambda)=\{a-a':a\in A,\; a'\in A',\; a>a'\},\]
and
\[ H(\widetilde{\lambda})=\{a'-a:a'\in A',\; a\in A,\; a'>a\}.\]
\end{enumerate}
\end{lem}

\section{The Splintering Lemma}\label{splintering}

The Splintering Lemma allows us to match up some of the entries of clustered partitions with complementary front sections easily. The key ingredient is the idea of \emph{$\infty$-partitions}; these $\infty$-partitions are not partitions in the usual sense, because they extend infinitely far to the left. They have a finite number, say $p$, of rows, the \emph{height} of the $\infty$-partition.
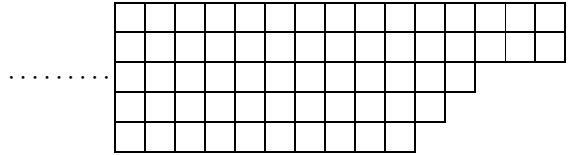
\[ \cdots\cdots\cdots\yng(15,15,12,11,10)\]
Notice that we can legitimately define the hook number multiset of this $\infty$-partition, since each number appears at most $p$ times. Given a partition $\lambda$, we can construct an $\infty$-partition $\lambda_\infty$ from it, by extending the partition infinitely far to the left. Notice that $\infty$-partitions are in one-to-one correspondence with all possible front sections in the obvious way, by removing columns with $p$ boxes in them except for the right-most column. In the example above, the corresponding front section is $(6,6,3,2,1)$. We mimic the definition of complementary front sections and say that two $\infty$-partitions are \emph{complementary} if their corresponding front sections are complementary. 

Consider a partition datum, say $[\lambda^{(f)},\lambda^{(r)},e_\lambda]_p$. Since there is a one-to-one correspondence between the set of all front sections of $p$ rows and the set of all $\infty$-partitions of height $p$, we can replace the front section $\lambda^{(f)}$ by the corresponding $\infty$-partition $\lambda_\infty$. Thus if $\lambda_\infty$ is the $\infty$-partition given by the diagram
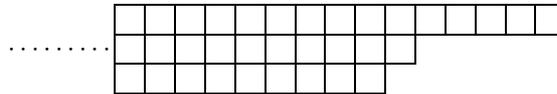
\[ \cdots\cdots\cdots\yng(15,10,9)\]
the partition with partition datum $[(7,2,1),(2),1]_3$ can be written as $[\lambda_\infty,(2),1]_3$. This is mererly a shorthand, since we will often consider $\infty$-partitions, and do not want to continually refer to the corresponding front section.

Write $H(\lambda_\infty)$ for the (infinite) multiset of hook numbers of the $\infty$-partition $\lambda_\infty$. The fundamental observation here is that after finitely many integers, each hook number appears exactly $p$ times in $H(\lambda_\infty)$, so we only need understand the `first part' of $H(\lambda_\infty)$; we could also think about $H(\lambda_\infty)$ in terms of the integers $\bar H(\lambda_\infty)$ missing from it (i.e., the multiset such that the union of it and the multiset of hook numbers together provide $p$ copies of every positive integer). For example, in the $\infty$-partition
\[ \lambda_\infty=\;\;\cdots\cdots\young(987521,6542,5431,321),\]
$\bar H(\lambda_\infty)=\{1,2,3,3,4,6\}$. We will call the multiset $\bar H(\lambda_\infty)$ the multiset of \emph{missing hook numbers} for the $\infty$-partition $\lambda_\infty$. We call two $\infty$-partitions $\lambda_\infty$ and $\mu_\infty$ of the same height \emph{clustered} if $H(\lambda_\infty)=H(\mu_\infty)$, or equivalently $\bar H(\lambda_\infty)=\bar H(\mu_\infty)$.

To describe an $\infty$-partition of $p$ rows, we give the $p$ hook numbers coming from the right-most column that has $p$ boxes in it a special name, the \emph{characteristic} of the $\infty$-partition. This is the same as the first-column hook numbers of the corresponding front section. In the example above, the characteristic is $(7,4,3,1)$. Note that the last element of the characteristic is always $1$.

Now we determine the missing hook numbers in terms of the characteristic of an $\infty$-partition. The first part of this result was essentially obtained by Frame, Robinson and Thrall in \cite{framerobthrall}, couched in the language of ordinary partitions.

\begin{lem}\label{missinghookno} Let $\lambda_\infty$ be an $\infty$-partition, with characteristic $(a_1,\dots,a_p)$. Then the missing hook numbers in row $i$ (counting from the top row) are
\[ \{ a_i-a_j:i<j\leq p\},\]
and consequently
\[ \bar H(\lambda_\infty)=\{ a_i-a_j: 1\leq i<j\leq p\}.\]
\end{lem}
\begin{pf} Let $(a_1,\dots,a_p)$ denote the characteristic, and write $A$ for the set of all $a_j$. The missing hook numbers in row $i$ are all less than $a_i$, so we can restrict our attention to $\lambda^{(f)}$. The hook numbers in row $i$ of $\lambda^{(f)}$ are, by Lemma \ref{oldresultonhooks},
\[ \{a_i-a': a'\in A',\; a'<a_i\},\]
where $A'$ is the complement of $A$ in the set $\{0,\dots,t-1\}$. The numbers between $1$ and $a_i$ that are not part of this set are therefore the numbers
\[ \{a_i-a: a\notin A',\; a<a_i\}=\{a_i-a: a\in A,\; a<a_i\}=\{a_i-a_j:i<j\leq p\}.\]
The second statement in the lemma now follows from this trivially.
\end{pf}

Notice that if two $\infty$-partitions $\lambda_\infty$ and $\mu_\infty$ are clustered, with characteristics $a_i$ and $b_i$, then the largest elements of $\bar H(\lambda_\infty)$ and $\bar H(\mu_\infty)$ are the same. But the largest elements of these multisets are $a_1-a_p$ and $b_1-b_p$ respectively. Since $a_p=b_p=1$, we must have $a_1=b_1$. This implies that, if $\lambda^{(f)}$ and $\mu^{(f)}$ correspond to clustered $\infty$-partitions, then they have the same number of boxes in the first row.

The Splintering Lemma will be an easy consequence of Lemma \ref{missinghookno}. Before we state it, we will describe a way to think of complementary $\infty$-partitions: consider a doubly-infinite version of a partition of height $p$, such as that pictured below.

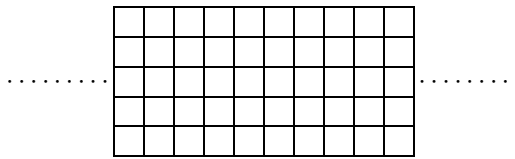
\[ \cdots\cdots\cdots \yng(10,10,10,10,10) \cdots\cdots\cdots \]

Now imagine snapping this like a piece of wood, so that the break turns the left-hand side into an $\infty$-partition and the right-hand side into an $\infty$-partition rotated by 180 degrees.
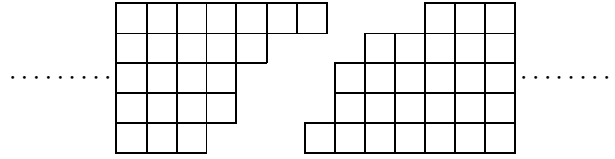
\[ \cdots\cdots\cdots \young(\,\,\,\,\,\,\,,\,\,\,\,\,,\,\,\,\,,\,\,\,\,,\,\,\,)\!\!\!\!\! \young(::::\,\,\,,::\,\,\,\,\,,:\,\,\,\,\,\,,:\,\,\,\,\,\,,\,\,\,\,\,\,\,) \cdots\cdots\cdots \]

These two $\infty$-partitions are complementary. It is an easy exercise for the reader to show that, if $(a_1,\dots,a_p)$ and $(b_1,\dots,b_p)$ denote the characteristics of two complementary $\infty$-partitions, then
\[ a_i=t-b_{p+1-i},\]
where $t$ is the row and column sum of the front section of either $\infty$-partition. (By a row and column sum, we mean the sum $t$ of the number of rows $r$ and the number of columns $c$, so that $t=r+c$.) We now give the Splintering Lemma.

\begin{thm}[Splintering Lemma] Let $\lambda_\infty$ and $\mu_\infty$ be complementary $\infty$-partitions, each with $p$ rows. Then $H(\lambda_\infty)=H(\mu_\infty)$.
\end{thm}
\begin{pf}This follows from Lemma \ref{missinghookno}: write $(a_1,\dots,a_p)$ for the characteristic of $\lambda_\infty$, and $(b_1,\dots,b_p)$ for the characteristic of $\mu$. Since $a_i=t-b_{p+1-i}$, we get
\begin{align*} \bar H(\mu_\infty)&=\{ b_i-b_j: 1\leq i<j\leq p\}
\\ &= \{(t-a_{p+1-i})-(t-a_{p+1-j}): 1\leq i<j\leq p\}
\\ &= \{a_{p+1-j}-a_{p+1-i}: 1\leq i<j\leq p\}
\\ &= \bar H(\lambda_\infty).\end{align*}
\end{pf}

There do exist non-complementary, clustered $\infty$-partitions; for example, the two $\infty$-partitions with front sections $(10,7,6,2,1^2)$ and $(10,8,7,6,1^2)$ have the same hook number multisets, but clearly are not complementary. The existence of non-complementary clustered front sections is crucial to the proof of Theorem \ref{mainthm2}.

\begin{prop}\label{twofisp} Let $\{\lambda,\mu\}$ be a period-$p$ cluster. Suppose that $\lambda$ and $\mu$ have remainders of the same weight, and the front sections, $\lambda^{(f)}$ and $\mu^{(f)}$, are complementary. Then $2|\lambda^{(f)}|\equiv 0\mod p$.
\end{prop}
\begin{pf} Since the remainders $\lambda^{(r)}$ and $\mu^{(r)}$ have the same weight, the rumps of $\lambda$ and $\mu$ must also have the same weight. Hence the front sections $\lambda^{(f)}$ and $\mu^{(f)}$ have the same weight modulo $p$, since in going from rumps to front sections, we remove a multiple of $p$ boxes. Hence
\[ |\lambda^{(f)}|\equiv |\mu^{(f)}|\mod p.\]
Finally, from the construction of complementary front sections,
\[ |\lambda^{(f)}|+|\mu^{(f)}|\equiv 0\mod p,\]
since they are constructed from a rectangle of height $p$. These two equations clearly imply that $2|\lambda^{(f)}|\equiv 0\mod p$, as required.
\end{pf}

\section{The Extension Lemma}\label{extension}

The Extension Lemma is a key result for our main theorem because it allows us to deduce the existence of periodic clusters from single clusters satisfying two natural conditions: if a cluster has clustered remainders and clustered $\infty$-partitions (for some period $p$), then the cluster is periodic. The periodic cluster given in the introduction with non-conjugate remainders had the weaker property of having clustered remainders, and the periodic cluster with non-complementary front sections had clustered $\infty$-partitions generated by those front sections; the author knows of no periodic clusters that do not satisfy both of these conditions.

\begin{lem}\label{sameextensions}
Suppose that $\lambda$ and $\mu$ are clustered partitions with partition data
\[ \lambda=[\lambda_\infty,\lambda^{(r)},e_\lambda]_p,\qquad \mu=[\mu_\infty,\mu^{(r)},e_\mu]_p.\]
In addition, suppose that $\lambda^{(r)}$ and $\mu^{(r)}$ are clustered, and that $\lambda_\infty$ and $\mu_\infty$ are clustered. Then $e_\lambda=e_\mu$.
\end{lem}
\begin{pf} Write $r_\lambda$ and $c_\lambda$, and $r_\mu$ and $c_\mu$, for the numbers of rows and columns of $\lambda$ and $\mu$ respectively. Then 
\[ t_\lambda=r_\lambda+c_\lambda=r_\mu+c_\mu=t_\mu.\]
Since the remainders $\lambda^{(r)}$ and $\mu^{(r)}$ are clustered, if we let $r_\lambda'$ denote the number of rows of $\lambda^{(r)}$, and so on, we have
\[ t_\lambda'=r_\lambda'+c_\lambda'=r_\mu'+c_\mu'=t_\mu'.\]
Next, $r_\lambda=r_\lambda'+p$, and $r_\mu=r_\mu'+p$. Now we need an expression for $c_\lambda$. Write $\lambda^{(f)}$ for the front section corresponding to $\lambda$, and similarly for $\mu^{(f)}$. Notice that the first row of both $\lambda^{(f)}$ and $\mu^{(f)}$ contain the same number of boxes, namely $c=t-p$, where $t$ is the row and column sum of the front section (see the discussion preceding the Splintering Lemma). We see that
\[ c_\lambda=c_\lambda'+e_\lambda+(c-1), \qquad c_\mu=c_\mu'+e_\mu+(c-1).\]
(It is $(c-1)$ rather than $c$ because the last row of any front section has one box in it, so the difference between the number of boxes in the first and last rows of a front section is $(c-1)$.)

We combine the expressions for $r_\lambda$ and $c_\lambda$, together with the equations $t_\lambda=t_\mu$ and $t_\lambda'=t_\mu'$ to get
\[ e_\lambda=e_\mu,\]
as claimed.
\end{pf}

\begin{thm}[Extension Lemma] Let $\bar\lambda$ and $\bar\mu$ be two clustered partitions, and let $p$ be an integer smaller than the number of rows both of $\bar\lambda$ and of $\bar\mu$. Let
\[ \bar\lambda=[\lambda_\infty,\lambda^{(r)},e_\lambda]_p,\qquad \bar\mu=[\mu_\infty,\mu^{(r)},e_\mu]_p\]
be the partition data. Suppose that $\lambda^{(r)}$ and $\mu^{(r)}$ are clustered, and that $\lambda_\infty$ and $\mu_\infty$ are clustered. Then, if $x$ is any non-negative integer, the two partitions
\[ \lambda=[\lambda_\infty,\lambda^{(r)},x]_p\text{ and }\mu=[\mu_\infty,\mu^{(r)},x]_p\]
are also clustered.
\end{thm}
\begin{pf}Firstly, we know that $e_\lambda=e_\mu$, by Lemma \ref{sameextensions}. Let $\lambda=[\lambda_\infty,\lambda^{(r)},x]_p$ and $\mu=[\mu_\infty,\mu^{(r)},x]_p$ be as in the statement, and let $\lambda'=[\lambda_\infty,\lambda^{(r)},0]_p$ and $\mu'=[\mu_\infty,\mu^{(r)},0]_p$. We will show that $H(\lambda)=H(\mu)$ if and only if $H(\lambda')=H(\mu')$. Then the fact that $H(\lambda)=H(\mu)$ when $x=e_\lambda$ will prove that it is true for all values of $x$.

Write $t'$ for the row and column sum of $\lambda'$, and $t$ for the row and column sum of $\lambda$: then
\[ t=t'+x.\]
Write $a_1,\dots,a_r$ and $x_1,\dots,x_r$ for the first-column hook numbers of $\lambda$ and $\lambda'$ respectively, so that, if $1\leq i\leq p$, then
\[ a_i=x_i+x,\]
and $a_i=x_i$ for $p+1\leq i\leq r$. Similarly, we write $b_1,\dots,b_s$ and $y_1,\dots,y_s$ for the first-column hook numbers of $\mu$ and $\mu'$ respectively, and we have similar relationships between the $b_i$ and the $y_i$.

Write $M_\lambda$ and $M_\mu$ for the multisets of hook numbers of $\lambda$ that lie in the rumps of $\lambda$ and $\mu$ respectively, and $M'_\lambda$ and $M'_\mu$ for the multisets of hook numbers of the rumps of $\lambda'$ and $\mu'$. We will have proven the result if we can show that $M_\lambda=M_\mu$ if and only if $M'_\lambda=M'_\mu$. Let $U_x$ denote the multiset consisting of $p$ copies of the integers between $1$ and $(t'-1)+x$ inclusive. (Then, for example, $U_0$ consists of $p$ copies of all integers between $1$ and $x_1$.) We have
\[ M_\lambda=M_\mu \iff U_x\setminus M_\lambda=U_x\setminus M_\mu.\]
We write $N_\lambda$ and $N_\mu$ for the complements of $M_\lambda$ and $M_\mu$ in $U_x$, and $N_\lambda'$ and $N_\mu'$ for the complements of $M_\lambda'$ and $M_\mu'$ in $U_0$. We therefore need to show that $N_\lambda=N_\mu$ if and only if $N_\lambda'=N_\mu'$. Since $M_\lambda$ is made up of the hook numbers of $\lambda$ that lie in the top $p$ rows, $N_\lambda$ is made up of the integers between $1$ and $a_1$ that are not in the $i$th row, for each $1\leq i\leq p$. This gives us a row-by-row decomposition of $N_\lambda$.

Let us use this decomposition of $N_\lambda$ to derive a description of it. Consider the $i$th row of $\lambda$. Then the integers between $1$ and $a_1$ that are not hook numbers lying in this row fall into two collections: those that lie between $1$ and $a_i$, and those that are between $a_i$ and $a_1$. The missing hook numbers between $1$ and $a_i$ are all $a_i-a_j$ for $i<j\leq r$, (by Lemma \ref{oldresultonhooks} or the discussion on $\infty$-partitions) and those above $a_i$ are simply all integers $\{a_i+1,\dots,a_1\}$. Write $R_i$ for the first multiset, and $S_i$ for the second. Then we have
\[ N_\lambda=\bigcup_{i=1}^p \( R_i\cup S_i\)=R_\lambda\cup S_\lambda,\]
where $R_\lambda$ and $S_\lambda$ have the obvious definition. Similarly, we construct $R_\mu$ and $S_\mu$. We can subdivide $R_\lambda$ into two disjoint multisets, by writing
\[ R_\lambda=\{a_i-a_j:1\leq i<j\leq p\} \cup \{a_i-a_j: 1\leq i\leq p,\,p+1\leq j\leq r\}.\]
Write $P_\lambda$ and $Q_\lambda$ respectively for the two multisets, and note that $P_\lambda=\bar H(\lambda_\infty)$. Similarly, write $R_\mu=P_\mu\cup Q_\mu$, and we get
\[ P_\lambda=\bar H(\lambda_\infty)=\bar H(\mu_\infty)=P_\mu.\]

Next, we perform the corresponding decomposition for $\lambda'$ and $\mu'$. Write
\[ N_\lambda'=P_\lambda'\cup Q_\lambda'\cup S_\lambda',\]
and $P_\mu'$, and so on for $\mu'$. Then
\[ P_\lambda'=P_\lambda=P_\mu=P_\mu'.\]
If we can show that $Q_\lambda\cup S_\lambda=Q_\mu\cup S_\mu$ if and only if $Q_\lambda'\cup S_\lambda'=Q_\mu'\cup S_\mu'$, then we will have shown that $M_\lambda=M_\mu$ if and only if $M_\lambda'=M_\mu'$, and we will have proven the theorem.

To this end, we will determine a relationship between $Q_\lambda$ and $Q'_\lambda$, and between $S_\lambda$ and $S'_\lambda$. The relationship between the last two multisets is easy to find, since
\begin{align*} S_\lambda&=\bigcup_{i=1}^p \{x+(x_i+1),\dots,x+(x_1)\}
\\ &= \{x+z: z\in S_\lambda'\}.\end{align*}
The relationship between $Q_\lambda$ and $Q_\lambda'$ is similar:
\begin{align*} Q_\lambda&=\{x+(x_i-a_j):1\leq i\leq p,\,p+1\leq j \leq r\}
\\ &=\left\{\bitbig x+z: z\in \{x_i-x_j:1\leq i\leq p,\,p+1\leq j \leq r\}\right\}
\\ &=\{x+z: z\in Q_\lambda'\}.\end{align*}
Thus
\[ Q_\lambda\cup S_\lambda=\{x+z: z\in Q_\lambda'\cup S_\lambda'\}.\]
A similar equation holds for $\mu$ and $\mu'$. Hence we clearly have
\[ Q_\lambda\cup S_\lambda=Q_\mu\cup S_\mu \iff Q_\lambda'\cup S_\lambda'=Q_\mu'\cup S_\mu',\]
which was what we wanted to prove. \end{pf}

The Extension Lemma allows us to prove that, for example, the cluster $\{A_n,A_n'\}$ given in Section 2 is actually a cluster, by proving that the $\infty$-parititons are clustered, the remainders are clustered, and that $A_0$ and $A_0'$ are clustered. Each of these is easy in this case, so it makes our goal of finding many periodic clusters much easier.

We now consider a partial converse to the Extension Lemma.

\begin{prop}\label{eitheror} Suppose that $\{\lambda,\mu\}$ is a periodic cluster, of period $p$. Let $\lambda^{(r)}$ and $\mu^{(r)}$ denote the remainders, and $\lambda_\infty$ and $\mu_\infty$ denote the $\infty$-partitions made from the rumps. Then $H(\lambda^{(r)})=H(\mu^{(r)})$ if and only if $H(\lambda_\infty)=H(\mu_\infty)$.
\end{prop}
\begin{pf}Recall that $H_i(\lambda)$ denotes the multiplicity of $i$ in the multiset $H(\lambda)$. Let $x$ be larger than any hook number in either $\lambda^{(r)}$ or $\mu^{(r)}$, and larger than any element of the missing hook numbers from $\lambda_\infty$ and $\mu_\infty$ (for example, let $x=|\lambda|$). Extend $\lambda$ and $\mu$ by $x$; so we may assume that, if $\lambda^{(t)}$ denotes the rump, then
\[ H_i(\lambda_\infty)=H_i(\lambda^{(t)})\]
for $i\leq x$, (and similarly for $\mu$). Since $\lambda$ and $\mu$ are clustered, for all $i\leq x$ we have $H_i(\lambda)=H_i(\mu)$.

Since we have extended the partitions, we now have
\[ H_i(\lambda)=H_i(\lambda_\infty)+H_i(\lambda^{(r)}).\]
Then we see that $H_i(\lambda_\infty)=H_i(\mu_\infty)$ if and only if $H_i(\lambda^{(r)})=H_i(\mu^{(r)})$. But the hypothesis is that one of these equations holds, and so they both do. This is true for all $i\leq x$, which is larger than both the entries in $\lambda^{(r)}$ and $\mu^{(r)}$ and the missing hook numbers of $\lambda_\infty$ and $\mu_\infty$, so if $\lambda_\infty$ and $\mu_\infty$ are clustered, then so are $\lambda^{(r)}$ and $\mu^{(r)}$, and vice versa.\end{pf}

This result explains why we get clustered remainders and clustered $\infty$-partitions appearing together. (Note that the result does not imply that either condition need hold: this is an open problem.)

\section{The Vertical Expansion Theorem}\label{vet}

In this section we prove a theorem that will substantially reduce the amount of work needed to prove Theorem \ref{mainthm2}. Specifically, if we can find a period-$p$ cluster of size $n$ and weight congruent to $d$ modulo $p$, then the Vertical Expansion Theorem asserts that we can find infinitely many different clusters with the same properties. Sections \ref{prepforthm} and \ref{proofofthm} are devoted to proving the existence of enough clusters to ensure that for each $n$, we can find a period $p$, and $p$ different period-$p$ clusters of size $2^n$, each with weight a different congruence class modulo $p$.

We start with a useful lemma, which lets us build up larger $\infty$-partitions from previously-known ones. Before we begin, let $\lambda_\infty$ and $\mu_\infty$ be two complementary $\infty$-partitions, with characteristics $(a_1,\dots,a_p)$ and $(t-a_p,\dots,t-a_1)$, where $t=a_1+1$. Then
\[ \{ a_i-a_j:1\leq i,j\leq p\}=\bar H(\lambda_\infty)\cup \{-x:x\in \bar H(\lambda_\infty)\}\cup p\cdot\{0\}.\]
To see this, $\bar H(\lambda_\infty)$ is exactly those $a_i-a_j$ such that $i<j$, the $p\cdot \{0\}$ comes from $a_i-a_j$ when $i=j$, and the $\{-x:x\in \bar H(\lambda_\infty)\}$ comes from the $-(a_j-a_i)$, recalling that
\[ \bar H(\lambda_\infty)=\bar H(\mu_\infty).\]
This will make our proof much easier.

\begin{lem}\label{twotimesinf} Let $\lambda_\infty$ and $\mu_\infty$ be two clustered $\infty$-partitions, with characteristics $(a_1,\dots,a_p)$ and $(b_1,\dots,b_p)$ respectively. Let $x$ be an integer at least as large as $a_1$. Write $\lambda_\infty'$ and $\mu_\infty'$ for the $\infty$-partitions with characteristics
\[ (a_1+x,\dots,a_p+x,a_1,\dots,a_p),\qquad \text{and}\qquad (b_1+x,\dots,b_p+x,b_1,\dots,b_p):\]
then $\lambda_\infty'$ and $\mu_\infty'$ are clustered.
\end{lem}
\begin{pf}Write $(a_1',\dots,a_{2p}')$ and $(b_1',\dots,b_{2p}')$ for the characteristics of $\lambda_\infty'$ and $\mu_\infty'$ respectively. (If $i>p$, then $a_i'=a_{i-p}$, and if $1\leq i\leq p$, then $a_i'=a_i+x$.) Then
\begin{align*} \bar H(\lambda_\infty')&=\{a_i'-a_j':1\leq i<j\leq 2p\}
\\ &=\{a_i'-a_j':1\leq i<j\leq p\}\cup\{a_i'-a_j':p+1\leq i<j\leq 2p\}\cup \{a_i'-a_j':1\leq i\leq p,\, p+1\leq j\leq 2p\}
\\ &=\{a_i-a_j:1\leq i<j\leq p\}\cup\{a_i-a_j:1\leq i<j\leq p\}\cup \{a_i-a_j+x:1\leq i,j\leq p\}.\end{align*}
The first two multisets are simply $\bar H(\lambda_\infty)$. The third multiset is that described in the discussion preceding this lemma, all of whose entries are incremented by $x$, and thus we clearly see that $\bar H(\lambda_\infty')$ is given by
\begin{align*} \bar H(\lambda_\infty')&=2\cdot \bar H(\lambda_\infty)\cup p\cdot\{x\} \cup \{x+a,x-a:a\in \bar H(\lambda_\infty)\}
\\ &=2\cdot \bar H(\mu_\infty)\cup p\cdot\{x\} \cup \{x+a,x-a:a\in \bar H(\mu_\infty)\}
\\ &=\bar H(\mu_\infty'),\end{align*}
as required.\end{pf}

This has the following technical corollary, which is a key step in the proof of the Vertical Expansion Theorem.

\begin{cor}\label{matchuptopbit} Let $\lambda_\infty$, $\mu_\infty$, $\lambda_\infty'$, $\mu_\infty'$, $a_i$, $b_i$, $a_i'$ and $b_i'$ be as in Lemma \ref{twotimesinf}. Then
\[ \{a_i'-a_j':1\leq i\leq p,\,i<j\leq 2p\}=\{b_i'-b_j':1\leq i\leq p,\,i<j\leq 2p\}.\]
\end{cor}
\begin{pf} The multiset is question is simply $\bar H(\lambda_\infty')$ without the multiset
\[ \{a_i'-a_j':p+1\leq i<j\leq 2p\}=\bar H(\lambda_\infty),\]
and clearly
\[ \bar H(\lambda_\infty')\setminus \bar H(\lambda_\infty)=\bar H(\mu_\infty')\setminus \bar H(\mu_\infty),\]
yielding the result.\end{pf}

\begin{prop}\label{vetpre} Let $\lambda=[\lambda_\infty,\lambda^{(r)},0]_p$ and $\mu=[\mu_\infty,\mu^{(r)},0]_p$ be two clustered partitions, with clustered remainders and clustered $\infty$-partitions. Let $\lambda'=[\lambda_\infty,\lambda,0]_p$ and $\mu'=[\mu_\infty,\mu,0]_p$. Then $\lambda'$ and $\mu'$ form a period-$p$ cluster.
\end{prop}
\begin{pf} The remainders of $\lambda'$ and $\mu'$ are clustered, as are the front sections of $\lambda'$ and $\mu'$, and so by the Extension Lemma, if $\lambda'$ and $\mu'$ are clustered then they form a periodic cluster. We will show that the rumps of $\lambda'$ and $\mu'$ have the same hook numbers. Write $(a_1',\dots,a_{r+p}')$ for the first-column hook numbers of $\lambda'$, and $(b_1',\dots,b_{s+p}')$ for the first-column hook numbers of $\mu'$.

Next, write $(a_1,\dots,a_r)$ for the first-column hook numbers of $\lambda$, and $(b_1,\dots,b_s)$ for the first-column hook numbers of $\mu$. Then $a_i$ and $a_i'$ are related by the equation
\[ a_i'=\begin{cases} a_{i-p}& i>p\\ a_i+(t-1)& i\leq p\end{cases},\]
where $t-1$ is the largest hook number in $\lambda^{(f)}$.

Using the same strategy as the proof of the Extension Lemma, we construct the multiset $U$, which consists of $p$ copies of every integer between 1 and $a_1'$ inclusive. Write $N_\lambda$ and $N_\mu$ for the complements of the hook numbers of the rumps of $\lambda'$ and $\mu'$ in the multiset $U$. We therefore need to show that $N_\lambda=N_\mu$. In a similar way to the proof of the Extension Lemma, we have an expression for $N_\lambda$, as
\[ N_\lambda=\{a_i'-a_j':1\leq i\leq p,\,i<j\leq p+r\}\cup\bigcup_{i=1}^p\{a_i'+1,\dots,a_1'\}.\]
The first multiset in this decomposition, say $R_\lambda$, is given by
\begin{align*} R_\lambda&=\{a_i'-a_j':1\leq i\leq p,\, i<j\leq 2p\} \cup \{a_i'-a_j':1\leq i\leq p,\,2p+1\leq j\leq p+r\}
\\ &=\{a_i'-a_j':1\leq i\leq p,\, i<j\leq 2p\} \cup \{a_i-a_j+(t-1):1\leq i\leq p,\,p+1\leq j\leq r\}.\end{align*}
By Corollary \ref{matchuptopbit}, the first multiset in this decomposition is equal to the corresponding multiset for $\mu'$, and so it remains to prove that
\begin{align*}\hphantom{=}&\{a_i-a_j+(t-1):1\leq i\leq p,\,p+1\leq j\leq r\}\cup\bigcup_{i=1}^p\{a_i+t,\dots,a_1+(t-1)\}
\\ =&\{b_i-b_j+(t-1):1\leq i\leq p,\,p+1\leq j\leq r\}\cup\bigcup_{i=1}^p\{b_i+t,\dots,b_1+(t-1)\}.\end{align*}

Now consider the period-$p$ cluster $\{\lambda'',\mu''\}$, where $\lambda''=[\lambda_\infty,\lambda^{(r)},t-1]_p$ and $\mu''=[\mu_\infty,\mu^{(r)},t-1]_p$. This time we know that the hook number multisets of the rumps of $\lambda''$ and $\mu''$ are the same. Notice that the first-column hook numbers of $\lambda''$ are equal to that of $\lambda$, except that the largest $p$ of them are incremented by $t-1$. (The same is true for $\mu''$.) Again, we take the complement of the hook numbers of the rumps in the suitable overset, consisting of $p$ copies of every integer between 1 and the largest hook number inclusive. We get the equation
\begin{align*} &\{\(a_i+(t-1)\)-\(a_j+(t-1)\): 1\leq i<j\leq p\}
\\ &\hphantom{=\;\:}\cup \{\(a_i+(t-1)\)-a_j: 1\leq i\leq p,\,p+1\leq j\leq r\}\cup \bigcup_{i=1}^p\{a_i+t,\dots,a_1+(t-1)\}
\\ &=\{\(b_i+(t-1)\)-\(b_j+(t-1)\): 1\leq i<j\leq p\}
\\ &\hphantom{=\;\:}\cup \{\(b_i+(t-1)\)-b_j: 1\leq i\leq p,\,p+1\leq j\leq r\}\cup \bigcup_{i=1}^p\{b_i+t,\dots,b_1+(t-1)\},\end{align*}
and since the first multisets in this equation are $\bar H(\lambda_\infty)$ and $\bar H(\mu_\infty)$, we get the exact formula that we needed to prove.\end{pf}

Notice that $\lambda'$ and $\mu'$ are clustered, and in fact are periodic of period both $p$ and $2p$. Thus there is no loss of generality in requiring the extensions of $\lambda$ and $\mu$ to be 0 in the statement of the proposition.

This proposition can be repeatedly applied to yield the Vertical Expansion Theorem.

\begin{thm}[Vertical Expansion Theorem] Suppose that $\lambda=[\lambda_\infty,\lambda^{(r)},e]_p$ and $\mu=[\mu_\infty,\mu^{(r)},e]_p$ are clustered partitions, such that the remainders and $\infty$-partitions are clustered. Let $(x_1,\dots,x_d)$ be a finite string of non-negative integers. Write $\lambda_0=\lambda^{(r)}$, and for each $1\leq i\leq d$, write
\[ \lambda_i=[\lambda_\infty,\lambda_{i-1},x_i]_p,\]
and similarly for $\mu_i$. Then $\{\lambda_d,\mu_d\}$ is a cluster, of period $jp$, for all $1\leq j\leq d$.
\end{thm}

This has an immediate corollary, offering the best-possible answer to a question of Herman and Chung in \cite{hermanchung}, namely whether one can find multiply-parameterized clusters. [This question was alluded to in the introduction.]

\begin{cor} There exist clusters parameterized by all finite strings of non-negative integers.\end{cor}

The best way to see the partition $\lambda_d$ in the statement of the Vertical Expansion Theorem is to have the partition $\lambda^{(r)}$ at the bottom, with $d$ copies of the front section of $\lambda_\infty$ bolted on top, shifted to the right so that the overhang from the $i$th one up to the $(i+1)$th is the quantity $x_i$.

To describe the partition $\lambda_d$ in the Vertical Expansion Theorem, we will extend the notation of partition data again, and write
\[ \lambda_d=[\lambda_\infty,\lambda^{(r)},(x_1,\dots,x_d)]_p.\]
The idea here is that successive copies of $\lambda^{(f)}$ are added on top of $\lambda^{(r)}$, each with extension $x_i$ for $1\leq i\leq d$. In particular, $\lambda=[\lambda_\infty,\lambda^{(r)},(e_\lambda)]_p$.

In order to find different clusters of the same weight, we need an algebraic description of the weight of a partition in terms of its partition datum.

\begin{lem}\label{weightofpartition} Suppose that $\lambda=[\lambda_\infty,\lambda^{(r)},e_\lambda]_p$ is a partition, and write $c^{(r)}$ for the number of columns of $\lambda^{(r)}$, and $\lambda^{(f)}$ for the front section corresponding to $\lambda_\infty$, as usual. Then
\[ |\lambda|=|\lambda^{(r)}|+|\lambda^{(f)}|+(c^{(r)}+e_\lambda-1)p.\]
\end{lem}

The proof of this is obvious, and left to the reader.

\begin{prop}\label{weightofbigpartition} Suppose that $\lambda=[\lambda_\infty,\lambda^{(r)},(x_1,\dots,x_d)]_p$ is a partition. Write $c^{(r)}$ for the number of columns of $\lambda^{(r)}$, and $\lambda^{(f)}$ for the front section corresponding to $\lambda_\infty$. Lastly, write $a$ for the first part of $\lambda^{(f)}$ minus 1. Then
\[ |\lambda|=|\lambda^{(r)}|+d|\lambda^{(f)}|+\left(dc^{(r)}+\frac{d(d-1)}2a+\sum_{i=1}^d(d+1-i)x_i-d\right)p.\]
In particular, modulo $p$,
\[ |\lambda|=|\lambda^{(r)}|+d|\lambda^{(f)}|.\]
\end{prop}
\begin{pf} Write $\lambda_0,\dots,\lambda_d$ as in the statement of the Vertical Expansion Theorem. To prove the proposition, we will count the number of boxes added in going from $\lambda_{i-1}$ to $\lambda_i$. This corresponds to the number of boxes in the top $p$ rows of $\lambda_i$. By Lemma \ref{weightofpartition}, since
\[ \lambda_i=[\lambda_\infty,\lambda_{i-1},x_i]_p,\]
we have
\[ |\lambda_i|=|\lambda_{i-1}|+|\lambda^{(f)}|+(b_{i-1}+x_i-1)p,\]
where $b_{i-1}$ is the number of columns in $\lambda_{i-1}$. We claim that for all $j>0$, this number is given by
\[ b_j=b_{j-1}+x_j+a,\]
and for $j=0$, it is given by $c^{(r)}$. Suppose that this is true. Then $b_1=c^{(r)}+x_1+a$, and in general
\[ b_j=c^{(r)}+(x_1+x_2+\cdots+x_j)+ja.\]
Hence we get
\[ |\lambda_i|=|\lambda_{i-1}|+|\lambda^{(f)}|+\left(\sum_{\alpha=1}^i x_\alpha+(i-1)a-1\right)p.\]

We can recursively apply this formula, noting that $|\lambda_0|=|\lambda^{(r)}|$, to get
\[ |\lambda|=|\lambda^{(r)}|+d|\lambda^{(f)}|+\left(dc^{(r)}+\frac{d(d-1)}2a+\sum_{i=1}^d(d+1-i)x_i-d\right)p.\]
It remains, therefore, to prove the assertion that we made on $b_j$. Certainly, $b_0=c^{(r)}$ since they are both defined to be the same thing, the number of columns in $\lambda^{(r)}$. For the inductive formula,
\[ b_j=b_{j-1}+x_j+a,\]
this is obviously true once we remember that $a$ is the number of columns in $\lambda^{(f)}$ minus 1.
\end{pf}

In particular, if $d=yp+1$ for some integer $y$, then $|\lambda_d|\equiv |\lambda|$ modulo $p$.

Suppose that $\Lambda$ is a period-$p$ cluster, and let $\lambda$ be a partition in $\Lambda$. Write $\lambda=[\lambda_\infty,\lambda^{(r)},e]_p$. From this, construct another partition, $\lambda'$, given by
\[ \lambda'=[\lambda_\infty,\lambda^{(r)},(\underbrace{1,1,\dots,1}_{p},e)]_p.\]
Let $\Lambda'$ be the period-$p$ cluster given by all such $\lambda'$, as $\lambda$ ranges over the partitions in $\Lambda$. Certainly $|\lambda|<|\lambda'|$, and by the remarks above, $|\lambda|\equiv |\lambda'| \mod p$. Thus there is an integer $x$ such that, if $\bar \lambda=[\lambda_\infty,\lambda^{(r)},x]_p$, then
\[ |\bar\lambda|=|\lambda'|.\]
We claim that $H(\bar \lambda)\neq H(\lambda')$. To see this, let us count the number of 1s occurring in the two multisets. Write $a$ for the number of occurrences of $1$ in the front section $\lambda^{(f)}$, and $b$ for the number of occurrences of $1$ in the remainder $\lambda^{(r)}$. It is easy to see that
\[ H_1(\bar\lambda)=a+b,\]
whereas
\[ H_1(\lambda')=(p+1)a+b,\]
proving the assertion, since $a\geq 1$.

Repeating this procedure, given one period-$p$ cluster of size $n$ and weight congruent to $d$ modulo $p$, we can find arbitrarily many clusters, each of size $n$ and weight congruent to $d$ modulo $n$, and each with different hook number multisets. However, Theorem \ref{mainthm2} stated that these clusters could be chosen so that not only are the hook numbers different for different clusters, but that the \emph{product} of those numbers is different. To show this, we need to use prime numbers, and period-$p^2$ clusters.

Let $\lambda=[\lambda_\infty,\lambda^{(r)},0]_p$ be a period-$p$ cluster, and consider the $p$ different period-$p^2$ clusters
\[ \lambda_i=[\lambda_\infty,\lambda^{(r)},(\underbrace{0,0,\dots,0}_{p},i)]_p,\]
as $i$ ranges between $1$ and $p$. It is easy to see using Proposition \ref{weightofbigpartition} that
\[ |\lambda_i|=p+|\lambda_{i-1}|.\]
Since all of the $|\lambda_i|$ are multiples of $p$, and cover all congruence classes modulo $p^2$ that are multiples of $p$ themselves, we must have
\[ |\lambda|\equiv |\lambda_h|\mod p^2,\]
for some $h$. [In fact, it is easy to show that $h=p$, but we do not need this.]

Let $x$ be an integer such that, if $\bar\lambda=[\lambda_\infty,\lambda^{(r)},x]_p$, then $|\bar\lambda|=|\lambda_h|$. Since $\bar\lambda$ is a period-$p$ cluster and $\lambda_h$ is a period-$p^2$ cluster, if we extend $\lambda_h$ by 1, we need to extend $\bar\lambda$ by $p$ in order to keep their weights the same. Write $t^{(r)}$ for the row and column sum of $\lambda^{(r)}$. Suppose that the largest hook number of $\lambda_h$ is $a$, and that the largest hook number of $\bar \lambda$ is $b$. Every time we extend $\lambda_h$ by 1, $a$ is increased by 1, whereas $b$ is increased by $p$, which for now we will assume is at least 3. [A similar procedure can deal with the uninteresting case $p=2$.] Write $\bar\lambda'$ and $\lambda_h'$ for the extended versions of these partitions, extended by $py$ and $y$ respectively. For all sufficiently large $y$, we have the inequality
\[ (b+py)-t^{(r)}\geq 2(a+y+1)+1.\]

Now we need Bertrand's Postulate, a famous result of Chebyshev, in 1852. This states that for any natural number $n$, between $n$ and $2n$ one can find a prime number. Applying this with $n=a+y+1$ yields a prime number $\ell$ such that 
\[ a'<\ell<b'-t^{(r)},\]
where $a'=a+y$ and $b'=b+py$. Since $a'<\ell$, certainly $\ell$ does not divide any of the elements of $H(\lambda_h')$, so in particular, does not divide their product. If we can show that it divides the product of the hook numbers of $\bar \lambda'$, we will have gone a significant way to proving our statement.

There are two ways to see that $\ell$ must divide one of the hook numbers. The first uses the so-called $\ell$-abacus, once we notice that the first-column hook numbers of $\bar \lambda'$ are either those of $\lambda^{(r)}$ (and hence at most $t^{(r)}$ or close to $b'$, and so there must be a space underneath the bead corresponding to $b'$ in the $\ell$-abacus. Since defining the abacus and explaining the concepts would be too complex, we provide another easy proof.

Consider the top row of $\bar \lambda'$, with largest hook number $b'$. The boxes directly above the remainder do not have consecutive hook numbers, but once we move along the first row until we no longer lie above the remainder, the hook numbers become consecutive all the way until the front section. By construction, the prime $\ell$ must fall in the region of consecutive hook numbers, since we made sure that $\ell<b'-t^{(r)}$. Hence $\ell$ is one of the hook numbers of $\bar\lambda'$, so divides their product.

In conclusion, we have shown that for all sufficiently large (period-$p^2$) extensions of $\bar\lambda$ and $\lambda_h$, we have that the products of their hook numbers are different. By extending the original partition $\lambda$ by $i$, where $1\leq i\leq p-1$, we get a period-$p$ partition with weight congruent to the other $\lambda_i$ modulo $p^2$, so we can employ the same procedure to get the following theorem.

\begin{thm}\label{differentproducts} Suppose that $\Lambda$ is a period-$p$ cluster, with weight congruent to $d$ modulo $p$. Then there exists, for all sufficiently large $n$ congruent to $d$ modulo $p$, two clusters, $\Lambda_1$ and $\Lambda_2$, one an extension of $\Lambda$, such that if $\lambda_1\in \Lambda_1$ and $\lambda_2\in \Lambda_2$, then $\lambda_1$ and $\lambda_2$ have different products of hook numbers.
\end{thm}

The construction used to find the two clusters above can be repeated, to yield the following corollary.

\begin{cor}\label{whatwewant} Let $\Lambda$ be a period-$p$ cluster of $n$ partitions, and weight congruent to $d$ modulo $p$. Then there exist arbitrarily many period-$p$ clusters $\Lambda_1,\Lambda_2,\dots$ all of size $n$ and of the same weight, and this weight is congruent to $d$ modulo $p$. Moreover, the $\Lambda_i$ can be chosen so that if $\lambda_i\in \Lambda_i$ and $\lambda_j\in \Lambda_j$, with $i\neq j$, then $\lambda_i$ and $\lambda_j$ have different products of hook numbers.
\end{cor}

\section{Preparing for the Proof: A New $\infty$-Partition}\label{prepforthm}

The Splintering Lemma implies the following result, which is essential for the proof of our main theorem.

\begin{lem}\label{comphooknos} Let $\lambda$ be a partition, and let $\widetilde{\lambda}$ denote its remnant. Write $t=r+c$, as usual. Then
\[ H_i(\lambda)-H_{t-i}(\lambda)=H_i(\widetilde{\lambda})-H_{t-i}(\widetilde{\lambda}).\]
\end{lem}
\begin{pf} Form the $\infty$-partition $\lambda_\infty$ of height $r$ by extending $\lambda$ infinitely far to the left. Let $\lambda'_\infty$ denote the complementary $\infty$-partition. Then $H(\lambda_\infty)=H(\lambda_\infty')$; in particular, $H_i(\lambda_\infty)=H_i(\lambda_\infty')$. Now we simply work out how $H_i(\lambda)$ and $H_i(\lambda_\infty)$ are related. Write $A$ for the set of first-column hook numbers of $\lambda$. Any occurrence of $i$ in $\lambda_\infty$ appears either in $\lambda$ itself, or to the left of $\lambda$, and so must be larger than the first-column hook number of $\lambda$ in that row. Thus
\[ H_i(\lambda)=H_i(\lambda_\infty)-|\{a\in A: a<i\}|,\]
and similarly
\[ H_{t-i}(\lambda)=H_{t-i}(\lambda_\infty)-|\{a\in A: a<(t-i)\}|.\]

On the other hand, consider $H_i(\lambda_\infty')$. Write $\lambda'$ for the front section of $\lambda_\infty'$, and notice that if we remove the left-hand column of $\lambda'$, we get the partition $\widetilde{\lambda}$. Let $B$ denote the first-column hook numbers of $\lambda'$; then
\[ B=\{t-a:a\in A\}.\]
Also, we have
\[ H_i(\widetilde{\lambda})=H_i(\lambda_\infty')-|\{b\in B:b\leq i\}|,\]
and
\[ H_{t-i}(\widetilde{\lambda})=H_{t-i}(\lambda_\infty')-|\{b\in B:b\leq(t-i)\}|.\]
The slight difference in the formulae comes from the fact that we do not want to include the hook numbers that make up $B$ when calculating $H_i(\widetilde{\lambda})$.

The final stage in the proof is to notice that, since $B=\{t-a:a\in A\}$, we have
\[ |\{b\in B:b\leq i\}|=|\{a\in A:a\geq (t-i)\}|,\qquad|\{b\in B:b\leq (t-i)\}|=|\{a\in A:a\geq i\}|.\]
Lastly, we see that
\[ |\{a\in A:a\geq (t-i)\}|=|A|-|\{a\in A:a<i\}|,\qquad |\{a\in A:a\geq i\}|=|A|-|\{a\in A:a<(t-i)|,\]
and upon collation of these facts, the result follows.\end{pf}

This result may be somewhat surprising: what it essentially says is that $H_i(\widetilde{\lambda})-H_{t-i}(\widetilde{\lambda})$ is determined by $H(\lambda)$, and so is the same for any clustered partitions. However, in \cite{hermanchung}, Herman and Chung show that $H(\lambda)$ and $H(\widetilde{\lambda})$ together determine the partition $\lambda$ up to conjugation. 

Lemma \ref{comphooknos} can be rewritten as
\[ H(\lambda)\setminus \{ t-h:h\in H(\lambda)\}=H(\widetilde{\lambda})\setminus\{t-h:h\in H(\widetilde{\lambda})\},\]
where the multisets involved can have negative multiplicities. This is easy to see since the multiplicity of $i$ in $\{ t-h:h\in H(\lambda)\}$ is equal to the multiplicity of $t-i$ in $H(\lambda)$, and similarly for $\widetilde{\lambda}$: hence the equation above, for the multiplicity of each $i$, becomes the equation given in Lemma \ref{comphooknos}.

Let $\lambda$ denote an arbitrary partition; we will construct a new $\infty$-partition $E_\infty(\lambda)$ from $\lambda$, of height $t$, where $t$ is the row and column sum. Let $A$ denote the set of first-column hook numbers of $\lambda$, and $B$ denote the set of top-row hook numbers of $\lambda$. Then the first-column hook numbers of the front section of $E_\infty(\lambda)$ are
Let $A$ denote the set of first-column hook numbers of $\lambda$, and $B$ denote the set of top-row hook numbers of $\lambda$. Then the first-column hook numbers of the front section of $E_\infty(\lambda)$ are
\[ \{a+t+1:a\in A\} \cup \{t-b:b\in B\}.\]
By Lemma \ref{oldresultonhooks}, this is the same as the set
\[ \{a+t+1:a \in A\} \cup \{a'+1:a'\in A'\},\]
where $A'$ is the complement of $A$ in the set $\{0,1,\dots,t-1\}$.

This $\infty$-partition is best-constructed by example. Let $t=r+c$, as usual, and construct the partition $\left((c+1)^t\right)$, a rectangle of height $t$ and width $c+1$. Remove from the bottom-right of the rectangle, the partition $\lambda$ reflected in the bottom-left to top-right diagonal. Then adjoin $\lambda$ to the top-right of this partition. This becomes the front section of the $\infty$-partition $E_\infty(\lambda)$. The example $\lambda=(4,2,1)$ is constructed below.
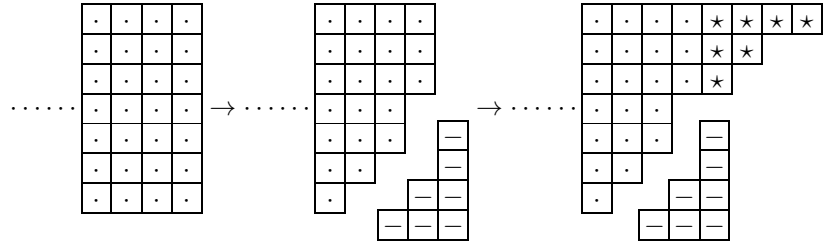
\[ \cdots\cdots\young(\cdot\cdot\cdot\cdot,\cdot\cdot\cdot\cdot,\cdot\cdot\cdot\cdot,\cdot\cdot\cdot\cdot,\cdot\cdot\cdot\cdot,\cdot\cdot\cdot\cdot,\cdot\cdot\cdot\cdot) \to \cdots\cdots\young(\cdot\cdot\cdot\cdot,\cdot\cdot\cdot\cdot,\cdot\cdot\cdot\cdot,\cdot\cdot\cdot,\cdot\cdot\cdot,\cdot\cdot,\cdot)\hspace{-0.3in}\raisebox{-0.65in}{\young(::-,::-,:--,---)} \to \cdots\cdots\young(\cdot\cdot\cdot\cdot\star\star\star\star,\cdot\cdot\cdot\cdot\star\star,\cdot\cdot\cdot\cdot\star,\cdot\cdot\cdot,\cdot\cdot\cdot,\cdot\cdot,\cdot)\hspace{-0.95in}\raisebox{-0.65in}{\young(::-,::-,:--,---)}\]
Here, a box with $\cdot$ in it is one that remains from the original rectangle, $-$ indicates that this box is removed, and $\star$ indicates that this is a box added. The boxes with $\star$ in them comprise a copy of the original partition $\lambda$.

\begin{prop}\label{infpartsofenvpart} Let $\lambda$ and $\mu$ be clustered partitions. Then $E_\infty(\lambda)$ and $E_\infty(\mu)$ are clustered $\infty$-partitions.
\end{prop}
\begin{pf} We will show that $\bar H(E_\infty(\lambda))$ is determined by $H(\lambda)$: then it must be true that
\[ \bar H(E_\infty(\lambda))=\bar H(E_\infty(\mu)).\]
Let $A$ denote the first-column hook numbers of $\lambda$. Then the first-column hook numbers of $E_\infty(\lambda)$ are
\[ C=\{a+t+1:a \in A\} \cup \{a'+1:a'\in A'\},\]
where again $A'$ denotes the complement of $A$ in the set $\{0,\dots,t-1\}$. Since $\bar H(E_\infty(\lambda))$ is the multiset of all differences between elements of $C$, we have
\begin{align*}\bar H(E_\infty(\lambda))&=\{ c_1-c_2:c_i\in C,\; c_1>c_2\}
\\ &=\{(a_1+t+1)-(a_2+t+1):a_i\in A,\; a_1>a_2\}\cup \{(a_1'+1)-(a_2'+1):a_i'\in A',\; a_1'>a_2'\}
\\ &\hphantom{=\;\:}\cup \{(a+t+1)-(a'+1):a\in A,\;a'\in A'\}
\\ &=\{a_1-a_2:a_i\in A,\; a_1>a_2\} \cup \{a_1'-a_2':a_i'\in A',\; a_1'>a_2'\} \cup \{a-a'+t:a\in A,\;a'\in A'\}
\\ &=\{a_1-a_2:a_i\in A,\; a_1>a_2\} \cup \{a_1'-a_2':a_i'\in A',\; a_1'>a_2'\}
\\ &\hphantom{=\;\:}\cup \{a-a'+t:a\in A,\;a'\in A',\; a>a'\}\cup \{a-a'+t:a\in A,\;a'\in A',\; a<a'\}.\end{align*}

The third multiset in this decomposition is $\{h+t:h\in H(\lambda)\}$, so is determined by $H(\lambda)$, and the fourth multiset in the decomposition is $\{t-h:h\in H(\widetilde{\lambda})\}$, where $\widetilde{\lambda}$ denotes the remnant of $\lambda$.

Now consider the multiset $X=\{ i-j: 0\leq j<i\leq t-1\}$. Clearly $X$ is dependent only on $H(\lambda)$, not on $\lambda$ itself. The multiset $X$ can be written as
\begin{align*} X=&\{a_1-a_2:a_i\in A,\; a_1>a_2\} \cup \{a_1'-a_2':a_i'\in A',\; a_1'>a_2'\}
\\ &\cup \{a-a':a\in A,\;a'\in A',\; a>a'\} \cup\{a'-a:a\in A,\;a'\in A',\; a'>a\}.\end{align*}
The third multiset in this decomposition is $H(\lambda)$, so that this is clearly only dependent on $H(\lambda)$. The fourth multiset in this decomposition is $H(\widetilde{\lambda})$. If we can show that the difference between the expressions for $\bar H(E_\infty(\lambda))$ and for $X$ is determined by $H(\lambda)$, then we are done. However, this is clearly true, since the difference between the two is
\[ \{t-h: h\in H(\widetilde{\lambda})\} \setminus H(\widetilde{\lambda}),\]
(allowing negative multiplicities of integers in this expression) and that this multiset is determined by $H(\lambda)$ is the expression following Lemma \ref{comphooknos}.\end{pf}

Notice that the multiset $Y=\{h+t:h\in H(\lambda)\}$ appears in $\bar H(E_\infty(\lambda))$. Indeed, by Lemma \ref{oldresultonhooks},
\[ H(\lambda)=\{a-a':a\in A,\;a'\in A',\; a>a'\},\]
and since the set of all elements of the characteristic of $E_\infty(\lambda)$ is
\[ \{a+t+1:a\in A\}\cup\{a'+1:a'\in A'\}=P\cup Q,\]
we clearly have $Y\subs \bar H(E_\infty(\lambda))$. Moreover, $Y$ is precisely those elements of $\bar H(E_\infty(\lambda))$ that are larger than $t-1$. To see this, write $P$ for the first multiset in the expression for the set of elements of the characteristic above, and $Q$ for the second. Let $x$ and $y$ denote two first-column hook numbers, and consider their difference, $x-y$. If $x$ and $y$ both lie in $P$, their difference is at most $t-2$. If $x$ and $y$ both lie in $Q$, their difference is likewise at most $t-2$. If $x$ lies in $P$ and $y$ lies in $Q$, then $x=a+t+1$ for some $a\in A$, and $y=a'+1$ for some $a'\in A'$. Now consider $x-y=t+(a-a')$: if $x-y\geq t$, then we must have $a\geq a'$, and since $a\neq a'$, we actually have $a>a'$, and so this is one of the elements of $Y$, as we asserted. Thus
\[ \{h+t:h\in H(\lambda)\}=\{h:h\in \bar H(E_\infty(\lambda)),\;h\geq t\}.\]

\section{The Enveloping Partition: The Proof of Theorem \ref{mainthm2}}\label{proofofthm}

Let $\lambda$ be a partition of $n$, and write $t=r+c$ for the sum of the rows and columns, as we have done previously. We denote by $E(\lambda)$ the partition with partition datum
\[ [E_\infty(\lambda),\lambda,0]_t.\]
This partition is called the \emph{enveloping partition} of $\lambda$. A better way to describe this partition is to take a $t\times t$ square, remove the reflection of $\lambda$ from the bottom-right, as in the construction of $E_\infty(\lambda)$, then add a copy of $\lambda$ to both the top-right corner and the bottom-left corner of the square. Thus, if $\lambda=(5,3,3,2)$ for example, $E(\lambda)$ is the partition

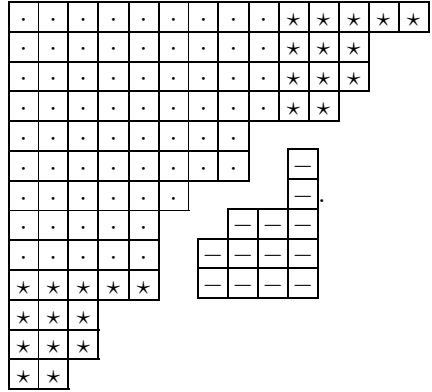
\[ \young(\cdot\cdot\cdot\cdot\cdot\cdot\cdot\cdot\cdot\star\star\star\star\star,\cdot\cdot\cdot\cdot\cdot\cdot\cdot\cdot\cdot\star\star\star,\cdot\cdot\cdot\cdot\cdot\cdot\cdot\cdot\cdot\star\star\star,\cdot\cdot\cdot\cdot\cdot\cdot\cdot\cdot\cdot\star\star,\cdot\cdot\cdot\cdot\cdot\cdot\cdot\cdot,\cdot\cdot\cdot\cdot\cdot\cdot\cdot\cdot,\cdot\cdot\cdot\cdot\cdot\cdot,\cdot\cdot\cdot\cdot\cdot,\cdot\cdot\cdot\cdot\cdot,\star\star\star\star\star,\star\star\star,\star\star\star,\star\star)\hspace{-1.2in}\raisebox{-0.5in}{\young(:::-,:::-,:---,----,----)}.\]
(Here, the $\cdot$ represents the boxes of the original square that remain, $-$ represents the boxes of the original square that are removed, and $\star$ represents the added copies of $\lambda$.)

It is clear from this picture that $E(\lambda^{(c)})=E(\lambda)^{(c)}$. Thus, if $\lambda$ and $\mu$ are conjugate, then $E(\lambda)$ and $E(\mu)$ are clustered. We can do much better.

\begin{thm}\label{mainthm3} Let $\lambda$ and $\mu$ be two clustered partitions. Then $E(\lambda)$ and $E(\mu)$ are clustered partitions.\end{thm}

We defer the proof of this result, but firstly deduce Theorem \ref{mainthm2} from it. Let $\lambda$ and $\mu$ denote two clustered partitions, possibly not periodic. Firstly, note that $\{E(\lambda),E(\mu)\}$ is a period-$t$ cluster: the remainders of $E(\lambda)$ and $E(\mu)$ are simply $\lambda$ and $\mu$, and so are clustered, and we have proven that $E_\infty(\lambda)$ and $E_\infty(\mu)$, the $\infty$-partitions of $E(\lambda)$ and $E(\mu)$, are clustered. Thus, by the Extension Lemma, $\{E(\lambda),E(\mu)\}$ is a period-$t$ cluster.

This enables us, given a cluster $\Lambda$ of $m$ partitions, to construct a cluster $E(\Lambda)_1$ of $2m$ partitions, by taking, for each $\lambda\in \Lambda$, the partition $E(\lambda)_1=[E_\infty(\lambda),\lambda,1]_t$, and the conjugate $E(\lambda)_1^{(c)}$ of $E(\lambda)_1$. The set
\[ E(\Lambda)_1=\{E(\lambda)_1,\,E(\lambda)_1^{(c)}:\lambda\in \Lambda\}\]
is then a period-$t$ cluster.

The next stage in the proof of Theorem \ref{mainthm2} is to notice that if $\lambda$ has weight $n$ and row and column sum $t$, then $E(\lambda)_1$ has weight $t^2+t+n$, period $t$, and row and column sum $3t+1$. We have an iterative procedure: if $\Lambda_1$ is simply a set consisting of a partition $\lambda$ and its conjugate, write $n_1$ for the weight of the partitions in $\Lambda_1$, and $t_1$ for the sum of the number of rows and the number of columns. Given a cluster $\Lambda_i$ of $2^i$ partitions, each with weight $n_i$ and row and column sum $t_i$, we construct a cluster $\Lambda_{i+1}$ of $2^{i+1}$ partitions, with weight
\[ n_{i+1}=n_i+t_i^2+t_i,\]
period $t_i$, and row and column sum
\[ t_{i+1}=3t_i+1,\]
by, for each $\lambda\in \Lambda_i$, considering the partition $E(\lambda)_1$, together with its conjugate. This means that, given any $\ell\in \N$, we can find a periodic cluster consisting of $2^\ell$ partitions. We need to show that we can find enough periodic clusters, each with the same period $t_{\ell-1}$, so that their weights cover all congruence classes modulo $t_{\ell-1}$. This would imply that for all sufficiently large $n$, there is a cluster with $2^\ell$ partitions in it, namely the periodic cluster with weight congruent to $n$ modulo $t_{\ell-1}$. By the Vertical Expansion Theorem (or rather, Corollary \ref{whatwewant}), we get arbitrarily many period-$t_{\ell-1}$ clusters with the correct weight modulo $t_{\ell-1}$. Then we will have Theorem \ref{mainthm2}.

To find enough clusters, firstly let $\ell$ be any integer, and write $t=t_1$, a variable. The period of the cluster $\Lambda_\ell$ above is $t_{\ell-1}=3t_{\ell-2}+1$, and is defined recursively. This gives
\[ t_{\ell-1}=3^{\ell-2}t+\frac{3^{\ell-2}-1}{2},\]
a linear function in $t$. Suppose that $t$ is odd, and write $t=r+c$, where $c=r-1$. The partitions with $r$ rows and $c$ columns have maximum weight $rc$, and minimum weight $r+c-1$, and for every integer between these two bounds, there is a partition with $r$ rows and $c$ columns of that weight. Notice also, that no partition with $r$ rows and $c$ columns is self-conjugate. We have
\[ r=\frac{t+1}{2},\qquad c=\frac{t-1}{2},\qquad rc=\frac{t^2-1}4,\qquad r+c-1=t.\]
Thus the difference $rc-(r+c-1)$ is given by
\[ \frac{t^2-4t-1}{4},\]
a quadratic function of $t$. Choose an odd $t$ such that
\[ \frac{t^2-4t-1}{4} > 3^{\ell-2}t+\frac{3^{\ell-2}-1}{2}=t_{\ell-1},\]
and let $\lambda^{(j)}$ denote a partition with $r$ rows, $c$ columns, and weight congruent to $j$ modulo $t_{\ell-1}$. (We know that such a partition exists by choice of $t$.) Finally, let
\[ \Lambda_1^{(j)}=\{\lambda^{(j)},{\lambda^{(j)}}^{(c)}\}.\]
Let $\Lambda_i^{(j)}$ denote the cluster obtained from $\Lambda_{i-1}^{(j)}$ in the way described above.

The clusters $\Lambda_i^{(j)}$ are each periodic, of period $t_{i-1}$ for all $0\leq j<t_{\ell-1}$ and all $2\leq i\leq \ell$. The weights of the clusters
\[ \Lambda_2^{(0)},\dots, \Lambda_2^{(t_{\ell-1}-1)}\]
cover all congruence classes modulo $t_{\ell-1}$, since to each weight, we have added $t_1^2+t_1$, and clearly by induction, since we add the same number to the weights of each cluster at each iteration, the weights of the clusters
\[ \Lambda_i^{(0)},\dots,\Lambda_i^{(t_{\ell-1}-1)}\]
cover all congruence classes modulo $t_{\ell-1}$. This implies that, if $N$ denotes the largest weight of the clusters $\Lambda_\ell^{(j)}$, then for all $n\geq N$, there is a period-$t_{\ell-1}$ cluster of size $2^\ell$ and weight $n$: Theorem \ref{mainthm2} follows.
\\ \qquad

It remains, therefore, to prove Theorem \ref{mainthm3}. This will be proven in a sequence of lemmas, which will show that the assertion that $H(E(\lambda))=H(E(\mu))$ (for a cluster $\{\lambda,\mu\}$) follows from the assertion that $E_\infty(\lambda)$ and $E_\infty(\mu)$ are clustered, a result that we already know. We provide an illuminating diagram for the proof of this theorem.

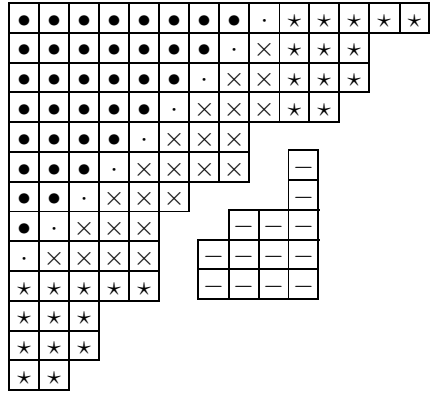
\[ \young(\bullet\bullet\bullet\bullet\bullet\bullet\bullet\bullet\cdot\star\star\star\star\star,\bullet\bullet\bullet\bullet\bullet\bullet\bullet\cdot\times\star\star\star,\bullet\bullet\bullet\bullet\bullet\bullet\cdot\times\times\star\star\star,\bullet\bullet\bullet\bullet\bullet\cdot\times\times\times\star\star,\bullet\bullet\bullet\bullet\cdot\times\times\times,\bullet\bullet\bullet\cdot\times\times\times\times,\bullet\bullet\cdot\times\times\times,\bullet\cdot\times\times\times,\cdot\times\times\times\times,\star\star\star\star\star,\star\star\star,\star\star\star,\star\star)\hspace{-1.2in}\raisebox{-0.5in}{\young(:::-,:::-,:---,----,----)}\]

The method of proof is the following: clearly, the two sets of boxes labelled with $\star$ are given by $2\cdot H(\lambda)$. The set of boxes labelled with $\cdot$ will be shown to have hook number $t$, and then we will show that the sum of the hook number in a box $(i,j)$ (with $i+j\leq t$) with $\bullet$, together with the hook number of the box $(t-j+1,t-i+1$)---the box reflected in the line made by the boxes with $\cdot$\ ---sum to $2t$. Finally, we show that the numbers denoted by $\bullet$ are simply the elements of $\bar H(E_\infty(\lambda))$, incremented by $t$. This will prove that $H(E(\lambda))$ is determined just by $H(\lambda)$, and does not require full knowledge of $\lambda$; thus $H(E(\lambda))=H(E(\mu))$ if $H(\lambda)=H(\mu)$.

For the remainder of the proof, let $\lambda$ denote a partition with $r$ rows and $c$ columns, write $t=r+c$, and let $E(\lambda)$ denote the enveloping partition. Write $h(i,j)$ for the hook number in the $(i,j)$ position of $E(\lambda)$. Write $A=\{a_1,\dots,a_r\}$ for the first-column hook numbers of $\lambda$, and again write $A'$ for the complement of $A$ in the set $\{0,\dots,t-1\}$.

\begin{lem}\label{firststep} The first-column hook numbers of $E(\lambda)$ are
\[ B=A\cup \{a'+t: a'\in A'\} \cup \{a+2t: a\in A\}.\]
\end{lem}
\begin{pf} Certainly $|B|=t+r$, which is the correct number. Thus we only have to show that the all of the elements in $B$ show up in the first-column hook numbers of $E(\lambda)$. This is safely left as an exercise for the reader, once we notice that $B$ is the union of $A$ and the first-column hook numbers of $E_\infty(\lambda)$, incremented by $(t-1)$.\end{pf}

\begin{lem}\label{secondstep} Let $1\leq i\leq t$. Then $h(i,t-i+1)=t$.\end{lem}
\begin{pf} Certainly, since $0\in A'$, we know that $t$ is a first-column hook number, by Lemma \ref{firststep}, and it is clearly $h(t,1)$. As $t-1\notin A'$, it must be that $2t-1$ is not a first-column hook number of $E(\lambda)$, and so $t$ is also a top-row hook number. Just as clearly, this top-row hook number must be $h(1,t)$. That the rest of the $h(i,t-i+1)$ are equal to $t$ is intuitively obvious from the diagram above, and we leave the reader to formulate a formal proof.
\end{pf}

\begin{lem}\label{thirdstep} Let $1\leq i,j\leq t$. Then
\[ h(i,j)+h(t-j+1,t-i+1)=2t.\]
\end{lem}
\begin{pf} By Lemma \ref{oldresultonhooks}, $h(i,j)+h(1,1)=h(i,1)+h(1,j)$ for all $1\leq i,j\leq t$. Then
\begin{align*} h(i,j)+h(t-j+1,t-i+1)&=\(\bitbig h(i,1)+h(1,j)-h(1,1)\)+\(\bitbig h(t-j+1,1)+h(1,t-i+1)-h(1,1)\)
\\ &= \(\bitbig h(i,1)+h(1,t-i+1)-h(1,1)\)+\(\bitbig h(1,j)+h(t-j+1,1)-h(1,1)\)
\\ &=h(i,t-i+1)+h(t-j+1,j)=2t.\end{align*}
\end{pf}

We restrict our attention to the hook numbers in the triangle where $i+j\leq t$. Recall that $h(1,1)=3t-1$, since the row and column sum is $3t$. Then
\[ h(1,j)+h(t-j+1,1)=h(1,1)+h(t-j+1,j)=4t-1,\]
and since we know the hook numbers $h(i,1)$ for $1\leq i\leq t$, this gives us an expression for $h(i,j)$, where $i+j\leq t$. Indeed,
\begin{align*} h(i,j)&=h(i,1)+h(1,j)-(3t-1)
\\ &=h(i,1)+\(\bitbig (4t-1)-h(t-j+1,1)\)-(3t-1)
\\ &=h(i,1)-h(t-j+1,1)+t.\end{align*}
Since $i+j\leq t$ we get
\[ \{h(a,b):a+b\leq t\}=\{h(a,1)-h(t-b+1,1)+t:\,a+b\leq t\}=\{h(i,1)-h(j,1)+t: 1\leq i<j\leq t\}.\]
This rewriting yields the following lemma.

\begin{lem}\label{fourthstep} We have
\[ \bar H(E_\infty(\lambda))=\{ h(i,j)-t: i+j\leq t\}=\{h(i,1)-h(j,1):1\leq i<j\leq t\}.\]
\end{lem}
\begin{pf} To prove this, we have to note that, if $(c_1,\dots,c_t)$ denotes the characteristic of $E_\infty(\lambda)$, then
\[ h(i,1)=c_i+(t-1).\]
This is simply the observation given in the proof of Lemma \ref{firststep}.\end{pf}

Now we can write down $H(E(\lambda))$.

\begin{lem}\label{fifthstep} The multiset $H(E(\lambda))$ is given by
\[ H(E(\lambda))=2\cdot H(\lambda)\cup \{ t+h: h\in \bar H(E_\infty(\lambda))\}\cup \left\{ t-h: h\in \bar H(E_\infty(\lambda))\setminus \{t+h:h\in H(\lambda)\}\right\}\cup t\cdot \{t\}.\]
\end{lem}
\begin{pf}Write $X$ for the multiset on the right-hand side of the formula. Since $\bar H(E_\infty(\lambda))$ contains $t(t-1)/2$ elements, the total number of elements in $X$ is $t^2+n$, where $n=|H(\lambda)|$. Thus we simply have to show that each element of $X$ shows up in $H(E(\lambda))$.

Certainly the two copies of $H(\lambda)$ show up, as these are the boxes with $\star$ in them in the diagram. Similarly, Lemma \ref{fourthstep} showed that the second term in $X$ is those boxes with $\bullet$ in them. The fourth term is the boxes with $\cdot$ in them, by Lemma \ref{secondstep}. It remains to discuss the third term.

In the discussion following Proposition \ref{infpartsofenvpart}, we showed that the elements of $\bar H(E_\infty(\lambda))$ that are larger than $t$ are precisely the elements
\[ \{t+h:h\in H(\lambda)\}.\]
These are those boxes with $\bullet$ in them that reflect onto those with $-$ in them, and since those boxes do not form part of $E(\lambda)$, we must remove them from the third term. Hence the third term in the description of $X$ is the multiset of hook numbers in the boxes with $\times$.

We have therefore showed that $X\subs H(E(\lambda))$, and since they have the same cardinality, we get the result.\end{pf}

This now establishes Theorem \ref{mainthm3}, since the multiset $H(E(\lambda))$ is determined by $H(\lambda)$, and so, in particular, if $H(\lambda)=H(\mu)$, then
\[ H(E(\lambda))=H(E(\mu)),\]
completing the proof of Theorem \ref{mainthm3}, and hence of Theorem \ref{mainthm2}.

\section*{Acknowledgements}

At the start of this work, Matt Towers and I had several discussions on this topic; he first noticed the complementary nature of the front sections of the clusters given in \cite{hermanchung}, and the fact that the remainders are clustered. In addition, he also conjectured the Splintering Lemma and the Extension Lemma. I would also like to thank my supervisor, Michael Collins, for his input on this subject, and John Wilson, for his stylistic comments.

\end{document}